\documentclass[a4paper,twoside]{amsart}
\usepackage{preamben}

\newtheorem{theorem}{Theorem}[section]
\newtheorem{lemma}[theorem]{Lemma}

\newtheorem{proposition}[theorem]{Proposition}
\newtheorem{corollary}[theorem]{Corollary}

\theoremstyle{definition}

\theoremstyle{remark}
\newtheorem{remark}[theorem]{Remark}

\numberwithin{equation}{section}

\figs

\newcommand{\mon}{\mathcal{M}\mathrm{on}}

\newcommand{\ve}{\varepsilon}
\newcommand{\tht}{\theta}

\newcommand{\Ad}{\mathrm{Ad}}

\newcommand{\trg}{\mathbb{T}}

\newcommand{\const}{Const}
\newcommand{\poi}{\mathcal{P}}
\newcommand{\real}{\mathrm{Re}}
\newcommand{\dis}{b}
\newcommand{\gabi}{\Psi}
\newcommand{\gabif}{\psi}
\newcommand{\tlgabi}{\widetilde{\Psi}}
\newcommand{\tlgabif}{\widetilde{\psi}}
\newcommand{\xze}{p}
\newcommand{\polg}{\psi}
\newcommand{\ircv}{\beta}
\newcommand{\bz}{\kappa}

\begin{document}
\title[The counterexample to a multidimensional Hilbert's...]{A counterexample to a multidimensional version of the weakened Hilbert's 16-th problem}
\vskip 1.truecm \author{Marcin Bobie\'nski} 
\author{Henryk \.{Z}o\l\k{a}dek}
\address{Institute of Mathematics, Warsaw University, ul. Banacha 2, 02-097
Warsaw, Poland}
\thanks{This research was supported by the KBN Grant No 2 P03A 015 29}
\email{mbobi@mimuw.edu.pl}
\email{zoladek@mimuw.edu.pl}

\date{\today}

\begin{abstract}
In the weakened 16th Hilbert's Problem one asks for a bound of the number of limit cycles which appear after a polynomial perturbation of a planar polynomial Hamiltonian vector field. It is known that this number is finite for an individual vector field. In the multidimensional generalization of this problem one considers polynomial perturbation of a polynomial vector field with invariant plane supporting a Hamiltonian dynamics. We present an explicit example of such perturbation with infinite number of limit cycles which accumulate at some separatrix loop.
\end{abstract}
\maketitle

\section{The result}
\label{sec:result}

Yu. Il'yshenko \cite{il16} and J. Ecalle \cite{ec} proved that an individual planar polynomial vector field can have only finite number of limit cycles.

On the other hand multi-dimensional vector fields with chaotic dynamics have infinite number of periodic trajectories. The Lorentz system \cite{mimr} and the Duffing system \cite{guho} provide best known examples. In the chaotic systems the periodic orbits are usually encoded by periodic sequences in a suitable symbolic dynamical system. This encoding is proved using topological methods (like the Lefschetz-Coneley index or Smale's horseshoe). This means that:
\begin{enumerate}
\item The periods of the periodic trajectories tend to infinity in rather irregular way.
\item The 1-cycles represented by different periodic trajectories have different ``topology'' i.e. they are linked between themselves.
\end{enumerate}
In particular, these cycles do not form a continuous family (so called center).

In Main Theorem below we give an example of polynomial 4-dimensional differential system, with infinite number of periodic solutions $\gamma_1,\gamma_2,\ldots$ such that
\begin{itemize}
\item the periods of $\gamma_j$ grow monotonically with $j$;
\item the corresponding 1-cycles have the same ``topology''; they are concentric cycles on an embedded invariant 2-dimensional disc of class $C^1$;
\item the $\gamma_j$ are isolated (they are limit cycles).
\end{itemize}

To construct the example we begin with the Hamiltonian planar system
\begin{equation}
  \label{hamsys}
  \dot x=X_H=(H_{x_2},-H_{x_1}), \quad (x_1,x_2)\in\bbR_x^2,\qquad H=x_1^3-3x_1 - x_2^2 + 2
\end{equation}
and the 2-dimensional linear system
\begin{equation}
  \label{2dy}
  \dot y= a y,\qquad y=y_1+i y_2\in\bbC \equiv \bbR^2_y,
\end{equation}
where $a=-\rho+i\omega$. Later we put $\rho=\omega=\sqrt3$.

\begin{figure}[htbp]
  \centering
  \input{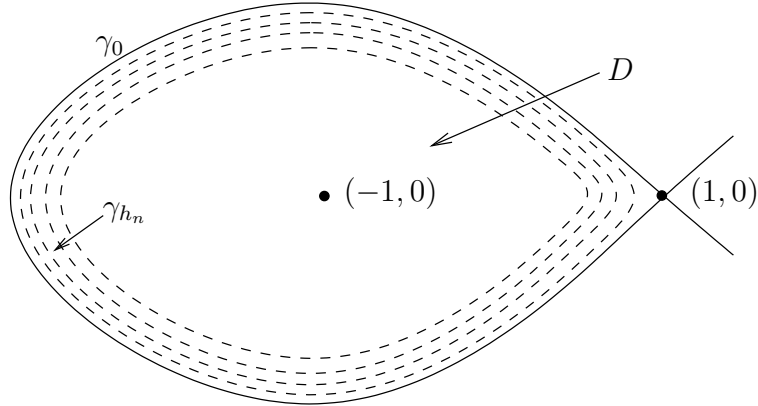}
  \caption{Phase portrait of the Hamiltonian vector field $X_H$ with ovals $\gamma_n$ generating limit cycles}
  \label{fig:ell}
\end{figure}

The Hamiltonian function from (\ref{hamsys}) is elliptic with the critical points $x=(-1,0)$ (center) and $x=(1,0)$ (saddle). The phase portrait of the field $X_H$ is shown on Figure \ref{fig:ell}

We consider the following coupling of the system (\ref{hamsys}) and (\ref{2dy})
\begin{equation}
  \label{system}
\left\{  \begin{aligned}
    \dot x &= X_H + \real(\overline{\bz}\,y)\, e_2\\
    \dot y &= a\, y + \ve H^4(x)\;(1-x_1),
  \end{aligned} \right.
\end{equation}
where $\ve>0$ is a small parameter, $e_2=(0,1)$ is a versor in $\bbR^2_x$ and $\bz\in\bbC$.

\begin{theorem}[Main]
\label{th:main}
Let $\rho=\omega=\sqrt3$ and
\begin{equation}
  \label{y0}
  \bz = 4\sqrt3 i + \frac{(3-3i)\sqrt6}{\sqrt\pi}(1+2i -\polg'(\tfrac{1-i}2)),
\end{equation}
where $\polg(z)$, the Euler Psi-function, is the logarithmic derivative of the Euler Gamma-function $\polg = \tfrac{\Gamma'}{\Gamma}$.

Then there exists an $\ve_0>0$ such that for any $0<\ve<\ve_0$ the system (\ref{system}) has a sequence of limit cycles $\gamma_n$, $n=1,2,\ldots$ which accumulate at the separatrix loop 
\begin{equation*}
\gamma_0=\{(x,y):\quad y=0,\ H(x)=0,\ x_1\leq1 \}
\end{equation*}
of the singular point $(x=(1,0),\,y=0)$ and lie on an invariant surface $y=\ve G(x,\ve)$ of class $C^1$.
\end{theorem}

\begin{remark}
The approximated numerical value of $\bz$ in formula (\ref{y0}) is
\[
\bz \approx -0.56 + 4.57 i.
\]
\end{remark}

Systems of the form
\begin{equation}
  \label{mdgensystem}
  \left\{\begin{aligned}
    \dot x &= X_H + F(x) y + \ve G(x)+\ldots \\
    \dot y &= A(x) y + \ve b(x) + \ldots,
  \end{aligned} \right.
\end{equation}
$x\in\bbR^2,\ y\in\bbR^\nu$, i.e. like (\ref{system}), appear in the so-called multidimensional generalization of the weakened 16-th Hilbert problem (see \cite{bo,bozoell,bozo3d,lezo}). Before perturbation, i.e. for $\ve=0$, we have the invariant plane $y=0$ with the Hamiltonian vector field $X_H(x)$. The ovals $H(x)=h$ form a 1-parameter family of its periodic trajectories. One asks how many of these trajectories survive the perturbation. In the 2-dimensional case ($\nu=0$) the linearization of the problem leads to the problem of real zeroes of an Abelian integral $I(h)=\int_{H(x)=h}\omega$; it is called the weakened 16-th Hilbert problem (see \cite{aril,il16}).

If $\nu\geq 1$, then the corresponding Pontryagin-Melnikov integrals (see \cite{mel,pont}), denoted $J(h)$, were found in \cite{lezo} and \cite{bozoell}. We call them the generalized Abelian integrals. 

The Abelian integrals $I(h)$ satisfy ODEs of the Fuchs type and have regular singularities with real spectrum (see \cite{ya}). Due to this, S. Yakovenko and others have found some effective estimations for the number of zeroes of $I(h)$. However, the generalized Abelian integrals do not satisfy any simple differential equation (see \cite{bo}) and sometimes have irregular singularities (e.g. at $h=\infty$). Moreover, even if the singularities are regular, then their spectra can be non-real.

Namely, the non-reality of the spectrum of $J(h)$ at the singularity $h=0$ is responsible for accumulation of zeroes of $J$. Below we find the asymptotics
\begin{equation*}
  J(h)  \sim C\, h^{9/2}\, \sin(\log\sqrt{h}),\qquad h\to 0^+.
\end{equation*}
It turns out that the zeroes $h_n\to 0^+$ of $J$ correspond to limit cycles $\gamma_n$ of the system (\ref{system}); the cycle $\gamma_n$ bifurcates form the oval $H^{-1}(h_n)$ (see Figure \ref{fig:ell}).

Therefore the system (\ref{system}) can be treated as a counterexample to the multi-dimensional weakened Hilbert's problem.

The remaining parts of the paper are devoted to the proof of Main Theorem. In Section \ref{sec:genab} we investigate the generalized Abelian integral and its zeroes. In Section \ref{sec:estim} we perform estimates needed for existence of genuine limit cycles.

\section{Proof of the Main Theorem}
\label{sec:proof}

\subsection{Generalized Abelian Integral}
\label{sec:genab}

The generalized Abelian integral is defined in two steps. Firstly one solves the so-called \emph{normal variation equation}
\begin{equation}
  \label{nveq}
  X_H(g) = a g + (1-x_1).
\end{equation}
Its solution $x\mapsto g(x)\in\bbC$ appears in the first (linear in $\ve$) approximation of the invariant surface (see the next section for more details)
\begin{equation*}
  y=\ve\,H^4(x)\; g(x) + O(\ve^2).
\end{equation*}

\begin{figure}[htbp]
  \centering
  \input{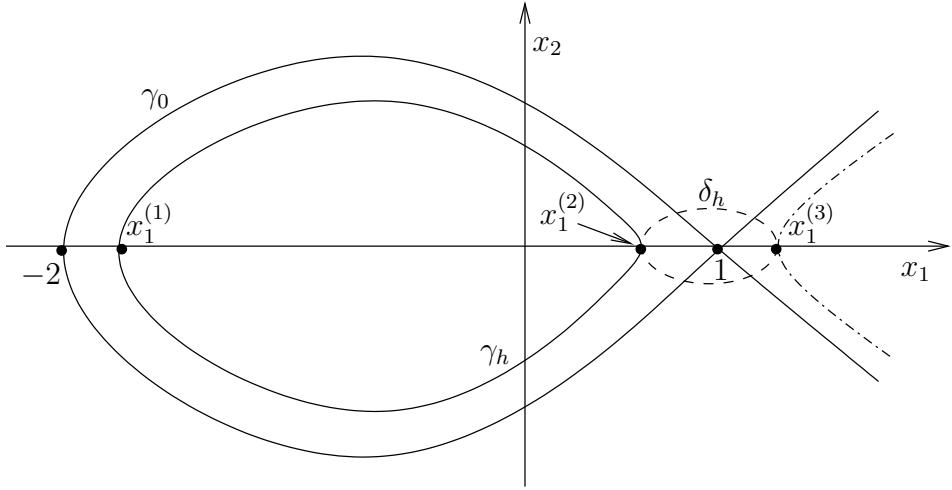}
  \caption{The basin $D$ filled with ovals $\gamma_h\subset H^{-1}(h)$.}
  \label{fig:dreg}
\end{figure}

We consider (\ref{nveq}) only in the basin $D=\{x:\ H(x)\geq 0,\, x_1\leq 1\}$ of the center $x=(-1,0)$, filled by the periodic solutions $\gamma_h (t)=\gamma (t)\subset \{H^{-1}(h)\}$, each of period
\begin{equation}
  \label{hamper}
  T_\gamma (h) = \int_{\gamma_h} -\frac{\dr x_1}{2 x_2} = \int_{\gamma_h} \dr t.
\end{equation}
We assume that the Hamiltonian time is chosen in such a way that for $0<h<4$ $x(0)=(x_1^{(1)},0)$, where $x_1^{(1)},x_1^{(2)},x_1^{(3)}$ are roots of the equation $H(x_1,0)-h=0$ (see Figure \ref{fig:dreg}). When restricted to $\gamma_h$, the equation (\ref{nveq}) is treated as the ODE $\dot g = a g + (1-x_1)$ with periodic boundary condition. Its unique solution is given in the integral form
\begin{equation}
  \label{gint}
  g(t,h) = (e^{-aT_\gamma} -1 )^{-1} \int_t^{t+T_\gamma} e^{a(t-s)} (1-x_1)(s,h)\; \dr s.
\end{equation}

Substituting the invariant surface equation (see Section \ref{sec:estim}) $y=\ve H^4 g + O(\ve^2)$ into the right hand side of $\dot x$ from (\ref{system}) we get the following perturbation of planar Hamiltonian system
\begin{equation}
  \label{hamres}
  \left\{\begin{aligned}
    \dot{x}_1 &= -2 x_2,\\
    \dot{x}_2 &= 3(1-x_1^2) + \ve\,H^4\,\real(\overline{\bz}\,g) + O(\ve^2).
  \end{aligned}\right.
\end{equation}
The generating function for limit cycles is given by the integral
\begin{equation}
  \label{point}
  J(h) = h^4 \int_{\gamma_h} \real\big(\overline{\bz}\,g(x)\big)\;\dr x_1.
\end{equation}
Let us denote the ``basic'' generalized Abelian integral (see \cite{bozoell}) by
\begin{multline}
  \label{psia}
  \gabi_\gamma (h) = \int_{\gamma_h} g(t) (1-x_1)(t) \dr t = \\
=(e^{-aT\gamma} -1 )^{-1} \int_0^{T_\gamma}\dr t \int_t^{t+T_\gamma}\dr s\; e^{a(t-s)} (1-x_1)(s)(1-x_1)(t).
\end{multline}
It is related to the generating function via the following
\begin{lemma}
  \label{lem:jviapsi}
  We have
\[
J(h) = h^4\, \real\Big[\overline{\bz}\,\Big(a\gabi_\gamma + 2\int_{\gamma_h} (1-x_1)\,\dr t \Big)\Big].
\]
\end{lemma}
\begin{proof}
In this proof we denote by dot, $\dot{f}=\tfrac{\dr}{\dr t} f = X_H(f)$, the differential with respect to the Hamiltonian time $t$. We have
\[
  J(h)= - h^4 \int_{\gamma} \real(\overline{\bz}\,g) \tfrac{\dr}{\dr t}{(1-x_1)}\, \dr t = h^4\, \real \Big(\overline{\bz}\int_{\gamma} \dot{g}(1-x_1)\, \dr t \Big).
\]
Next, $\dot{g}=a\,g+(1-x_1)$ gives
\begin{multline*}
 J(h)= h^4\, \real\Big[\overline{\bz}\,\Big( a\gabi_\gamma + \int_{\gamma} (1-x_1)^2\, \dr t \Big)\Big] =  h^4\, \real\Big[\overline{\bz}\Big( a\gabi_\gamma +  \int_{\gamma} (2-2 x_1 - \tfrac13 \dot{x}_2)\, \dr t \Big)\Big] =\\
  =h^4\, \real\Big[\overline{\bz}\Big(a\gabi_\gamma + 2\int_{\gamma_h} (1-x_1)\,\dr t \Big)\Big].
\end{multline*}
\end{proof}

Our next aim is to determine the leading terms in the asymptotic expansion as $h\to 0^+$ of the integrals $\int_{\gamma} (1-x_1)(t)\dr t$, $T_\gamma$ and $\gabi_\gamma$. We begin with the Abelian integrals. It is known \cite{zolbook} that these integrals extends to multivalued holomorphic functions with logarithmic singularities. We shall need explicit form of the leading terms.
\begin{lemma}
\label{lem:abi}
There exists an open neighborhood $0\in U\subset \bbC$ in the complex domain and holomorphic functions $\eta_0,\eta_1,\zeta_0,\zeta_1\in\Omega(U)$ such that
\begin{align}
  T_\gamma &= \eta_0(h) + \zeta_0(h)\, \log h = - \tfrac1{2\sqrt{3}} \log h+\tfrac{\sqrt3}2\log12 + O(h\log h),\label{tgex}\\
  \int_{\gamma} (1-x_1)\dr t &= \eta_1(h) + \zeta_1(h)\, \log h = 2\sqrt3 + O(h\log h)\label{abex}.
\end{align}
\end{lemma}
\begin{proof}
  We consider the pair of basis elliptic Abelian integrals
\begin{equation*}
  I_0(h) = T_\gamma = \int_{\gamma}\frac{-\dr x_1}{2 x_2},\qquad I_1(h) = \int_{\gamma}\frac{-x_1 \dr x_1}{2 x_2}.
\end{equation*}
Note that $\int_\gamma(1-x_1)\dr t = I_0-I_1$ and that $I_0=\tfrac{\dr}{\dr h}\Big(\text{area of } \{H>h\}\Big)$. 

These functions $(I_0,I_1)$ satisfy the Picard-Fuchs equations
  \begin{equation}
    \label{fuchs}
    \begin{split}
      6 h(h-4) I_0' &= -(h-2) I_0 - 2 I_1 \\
      6 h(h-4) I_1' &=  2 I_0  +(h-2) I_1. \\
    \end{split}
  \end{equation}
The other, independent solution to this system is the pair $(K_0,K_1)$, where
\begin{equation*}
  K_0(h) =  \int_{\delta_h}\frac{-\dr x_1}{2 x_2},\qquad K_1(h) = \int_{\delta_h}\frac{-x_1 \dr x_1}{2 x_2}
\end{equation*}
are integrals along another cycle $\delta_h$ in the complex curve $E_h=\{H(x)=h\}\subset\bbC^2$. If $h\in(0,4)$ then the polynomial $x_1^3-3 x_1 +2 -h$ has three real roots $x^{(1)}_h<x^{(2)}_h<x^{(3)}_h$ (see Figure \ref{fig:dreg}). The the cycle $\gamma_h$ (respectively $\delta_h$) is represented as the lift to the Riemann surface $E_h$ of loops in the complex $x_1$-plane surrounding the roots $x^{(1)}_h$ and $x^{(2)}_h$ (respectively $x^{(2)}_h$ and $x^{(3)}_h$). Note the following integral formulas for $T_{\gamma}(h)$:
\begin{equation}
  \label{tgi}
  T_{\gamma} = \int_{x^{(1)}_h}^{x^{(2)}_h} \frac{\dr x_1}{\sqrt{x_1^3-3 x_1 +2 -h}} = \int_{x^{(3)}_h}^\infty \frac{\dr x_1}{\sqrt{x_1^3-3 x_1 +2 -h}},\qquad h\in(0,4).
\end{equation}
The second equality corresponds to unobstructed deformation of integration contour $\gamma_h$ to loop surrounding $x^{(3)}_h$ and $\infty$.

The system (\ref{fuchs}) has resonant singular point $h=0$. Any its solution is either analytic near $h=0$ (it is $(K_1,K_2)$) or it has the form like $I_0,I_1$:
  \begin{equation}
    \label{i12exp}
    \begin{split}
     I_0(h) &= (a_0+ a_1 h + \ldots) +\tfrac1{2\pi i} K_0\;\log h,\\
     I_1(h) &= (b_0+ b_1 h + \ldots) +\tfrac1{2\pi i} K_1\;\log h.
    \end{split}
  \end{equation}
This representation follows from the Picard-Lefschetz formula
\begin{equation}
  \label{piclef}
  \gamma_h \longrightarrow \gamma_h\cdot\delta_h, \qquad\qquad\delta_h \longrightarrow \delta_h,
\end{equation}
which describes the monodromy transformations of the generators of $\pi_1(E_h,*)$, as $h$ surrounds the critical value $0$; here $*$ denotes a basepoint.

We need to calculate the expansions of $I_0,I_1$. As we shall see, it is enough to calculate $K_0(0)$ and $a_0$; all other coefficients follows from the system (\ref{fuchs}) and can be recursively determined. Indeed, to compensate terms with $\log h$ in (\ref{fuchs}) we must have
\begin{equation}
\label{auxk1k0}
K_1(0)=K_0(0).
\end{equation}
Terms with $h^0$ give
\begin{equation}
\label{auxb0a0}
b_0=a_0 +\tfrac{12}{2\pi i}K_0(0).
\end{equation}
It can be continued further.

To determine $K_0(0)$ and $a_0$ simultaneously, we make a coordinate change $u=(x_1-1)/(x^{(3)}_h-1)$ in the integral (\ref{tgi}); we denote also $\xze=(x^{(3)}_h-1)$. Since $\xze = \sqrt{h/3} +O(h)$ as $h\to 0^+$, the following integral 
\begin{multline*}
  \int_1^\infty \dr u\Big[\frac1{\sqrt{u^2(3+ \xze u) - h/\xze^2}} - \frac1{\sqrt{u^2(3+ \xze u)}} - \frac1{\sqrt{3 u^2-3}} + \frac1{\sqrt{3}} \Big] \xrightarrow{\ h\to 0^+} 0.
\end{multline*}
We calculate
\begin{align*}
  &\int_1^\infty \dr u\Big[ \frac1{\sqrt{3 u^2-3}} - \frac1{\sqrt{3}} \Big] = \frac{\log 2}{\sqrt 3}, \\
  &\int_1^\infty \dr u \Big[ \frac1{u \sqrt{3+ \xze u}} \Big] = \tfrac2{\sqrt 3}\log\left(\tfrac{2\sqrt3}{\sqrt{\xze}} + o(1)\right) = - \tfrac1{2\sqrt{3}} \log h + \frac{\log(12\sqrt3)}{\sqrt3} + o(h^{1/2}).
\end{align*}
Thus $a_0=\tfrac{\sqrt3}{2}\log 12$, $K_0(0)=- \tfrac{2\pi i}{2\sqrt3}$. Substituting these values to the relations (\ref{auxk1k0}), (\ref{auxb0a0}) and using the expansion (\ref{i12exp}) we get the leading terms of the expansions as in formulas (\ref{tgex}) and (\ref{abex}).\\
\end{proof}

Let us pass to expansion of $\gabi_\gamma$.
\begin{proposition}
\label{pr:genabi}
Let $-2\sqrt3 <\real(a)<0$. There exists an open neighborhood $0\in U\subset \bbC$ in the complex domain and holomorphic functions $\varphi_1,\varphi_2,\varphi_3$ such that
\begin{multline}
\label{genabian}
\gabi_\gamma (z) = \varphi_1(h) +\varphi_2(h)\, \log h + \varphi_3(h)\cdot \Big(e^{-a T_\gamma}-1\Big)^{-1} = C_0+C_1 h^{-a/2\sqrt3} + \ldots,
\end{multline}
where 
\begin{align}
  C_1 &= \frac{(\pi a)^2}{\sin^2(\pi a/2\sqrt3)}, \label{gabic1}\\
  C_0 &= \frac{3\sqrt2}{\sqrt \pi} \Big(-1+2w+2w^2 \polg'(-w)\Big), \qquad w=\tfrac{a}{2\sqrt3}. \label{gabic0}
\end{align}
\end{proposition}

\begin{remark}
\label{rk:y0viac0}
One can easily observe that the value (\ref{y0}) of $\bz$ satisfies the relation
\begin{equation}
  \label{y0viac0}
  \bz = i(4\sqrt3 + a C_0).
\end{equation}
It is chosen in a way to annihilate the leading term ($\sim h^4$) of $J(h)$ and to reveal the term with the infinite sequence of zeroes -- see the Corollary \ref{cor:zeroes} and its proof below.
\end{remark}

\begin{corollary}
\label{cor:zeroes}
Providing the values of parameters $(\rho,\omega,\bz)$ as in Theorem Main, the integral $J(h)$ (\ref{point}) has a sequence $h_n$, $n=1,2,\ldots$ of simple zeroes accumulating at $h=0$.
\end{corollary}
\begin{proof}
We calculate the leading term of the expansion of $J(h)$ using Lemma~\ref{lem:jviapsi}, Lemma~\ref{lem:abi} and Remark~\ref{rk:y0viac0}:
\begin{align*}
  J(h) &= h^4\real\Big[\overline{\bz}\,\big(a C_0 + a C_1 h^{1/2-i/2} + 4\sqrt3 +o(h^{3/4})\big)\Big] = \\
&=h^4\real\Big[\overline{\bz}\,\big(a C_0 + 4\sqrt3\big)\Big] + h^4\real\Big[\overline{\bz}\,a\, C_1 h^{1/2-i/2}\Big] + o(h^{4+3/4})=\\
&=R\, h^{4+1/2} \cos(\log \sqrt{h}-\alpha_0) + o(h^{4+3/4}),
\end{align*}
where $R=|\overline{\bz}\, a\, C_1|,\quad \alpha_0 = \mathrm{arg}(\overline{\bz}\, a\, C_1)$. Analogously we get
\[
J'(h) = R_1\, h^{3+1/2} \cos(\log \sqrt{h}-\alpha_1) + o(h^{3+3/4}).
\]
Thus, by Implicit Function Theorem the zeroes $(h_n)$ of $J(h)$ approximate the simple zeroes $h_n^{(0)}$ of the function $\cos(\log \sqrt{h}-\alpha_0)$.\\
\end{proof}

The remaining part of this section is devoted to the proof of Proposition \ref{pr:genabi}. It goes in two steps. In first one we show that the function 
\[
\gabi_\gamma(h) \longrightarrow C_0
\]
as $h\to 0^+$ and so is bounded. 

In the second step we determine the monodromy of the generalized Abelian integral $\gabi_\gamma(h)$ as $h$ surrounds $h=0$. We know that then $\gamma$ changes to $\mon_0\gamma=\gamma\cdot\delta$ and $\mon_0\delta=\delta$. We would like to express $\gabi_{\gamma\cdot\delta}$ in simple terms, in order to determine the singularity of $\gabi_\gamma(h)$ at $h=0$. Rather complicated formulas for $\gabi_{\gamma\cdot\delta}$ are given in \cite{bozoell} and \cite{bozo3d}. In \cite{bo} these formulas were simplified using certain upper triangle representation $\rho$ of the fundamental group $\pi_1(E_h,*)$. We recall this construction below.

We denote
\begin{equation}
\label{thtpdef}
  \begin{split}
    \gabi_\gamma (h) & = (e^{-aT\gamma} -1 )^{-1} \int_0^{T_\gamma}\dr t \int_t^{t+T_\gamma}\dr s\; e^{a(t-s)} (1-x_1)(s)\,(1-x_1)(t),\\
    \lambda_\gamma(h) & = e^{-a T_\gamma/2},\\
    \phi_\gamma(h) & = \lambda_\gamma \,\int_0^{T_\gamma}\dr t\int_0^t\dr s\; (1-x_1)(t)\, (1-x_1)(s)\cdot  e^{a(t-s)}= \\ 
    & =\lambda_\gamma \,\int_{-T_\gamma/2}^{T_\gamma/2}\dr t\int_{-T_\gamma/2}^t\dr s\; (1-x_1)(t+\tfrac{T_\gamma}2)\, (1-x_1)(s+\tfrac{T_\gamma}2)\,  e^{a(t-s)},\\
    \tht_\gamma^+ & =  \lambda_\gamma\,\int_0^{T_\gamma}\dr t\; (1-x_1)(t)\, e^{a t} = \int_{-T_\gamma/2}^{T_\gamma/2}\dr t\; (1-x_1)(t+\tfrac{T_\gamma}2)\, e^{a t},  \\
    \tht_\gamma^- & =  \lambda_\gamma^{-1}\,\int_0^{T_\gamma}\dr t\; (1-x_1)(t)\, e^{-a t} = \int_{-T_\gamma/2}^{T_\gamma/2}\dr t\; (1-x_1)(t+\tfrac{T_\gamma}2)\, e^{-a t},  
  \end{split}
\end{equation}
Here the subscript $\gamma$ underlines dependence of the above functions on the loop $\gamma=\gamma_h$.

We introduce the following space of triangular matrices
\begin{equation}
  \label{gform}
\trg \colon = \left\{ \left(\begin{smallmatrix}
    \lambda& \tht^{\scriptscriptstyle-} &\phi\\
    0& \lambda^{\scriptscriptstyle-1}& \tht^{\scriptscriptstyle+} \\
    0& 0& \lambda \\
    \end{smallmatrix}\right), \quad \lambda\in\bbC^*,\quad \tht^+,\tht^-,\phi\in\bbC \right\};
\end{equation}
it forms a group. For $W\in \trg$ we denote
\[
|W| = \det W = \lambda.
\]
Existence of a 2-dimensional Jordan cell is measured in the following formula
\begin{equation}
  \label{psirdef}
\frac{(W-|W|)(W-1/|W|)}{|W|^2-1} = \gabif(W)\left(\begin{smallmatrix}0&0&1\\ 0&0&0\\ 0&0&0\\ \end{smallmatrix}\right);
\end{equation}
explicitly we have
\begin{equation}
  \label{gabirexplicit}
  \gabif(W) = \frac{\tht^+\tht^-}{\lambda^2-1} + \frac{\phi}{\lambda}.
\end{equation}

\begin{theorem}[\cite{bo}]
\label{th:bo}
The map $\rho:\pi_1(E_h,*)\rightarrow \trg$,
\begin{equation}
  \label{thrho}
  \rho(\gamma) = \begin{pmatrix}
    \lambda_\gamma& \tht_\gamma^- & \phi_\gamma\\
    0& \lambda^{-1}_\gamma& \tht^+_\gamma \\
    0& 0& \lambda_\gamma \\
  \end{pmatrix},
\end{equation}
where $\lambda_\gamma,\tht_\gamma^\pm,\phi_\gamma$ are defined in (\ref{thtpdef}), defines a representation of the fundamental group of $E_h$. Moreover, we have
\begin{equation}
  \label{thgabi}
  \gabi_\gamma =\gabif\circ \rho\; (h).
\end{equation}
\end{theorem}
\begin{proof}[Sketch of the proof]
We have $\gabi_\gamma=(\lambda_\gamma^2-1)^{-1}\iint e^{a(t-s)} (1-x)(t)\;(1-x)(s)$, where the integration domain is $\Sigma=\{(t,s):\ 0\leq t \leq T_\gamma,\ t\leq s\leq T_\gamma+t\}$ (see (\ref{thtpdef})). We divide $\Sigma$ into two ``triangles'' $\triangle_1=\{0\leq t\leq T_\gamma,\ t\leq s\leq T_\gamma\}$ and $\triangle_2=\{0\leq t\leq T_\gamma,\ T_\gamma\leq s \leq T_\gamma+t\} = \triangle_0 + (0,T_\gamma)$, where $\triangle_0=\{0\leq t\leq T_\gamma,\ 0\leq s \leq t\}$.
We have $\iint_{\triangle_1+\triangle_2}(\cdot)= \iint_{\triangle_1+\triangle_0}(\cdot) + \iint_{\triangle_2-\triangle_0}(\cdot)$, where $\iint_{\triangle_1+\triangle_0}(\cdot)= \tht^+_\gamma\, \tht^-_\gamma$ and $\iint_{\triangle_2-\triangle_0}(\cdot)=(\lambda_\gamma^2-1)\lambda_\gamma^{-1}\phi_\gamma$. Now the formula (\ref{thgabi}) follows from (\ref{gabirexplicit}).

The property $\rho(\gamma\cdot\delta)=\rho(\gamma)\;\rho(\delta)$, $\gamma,\delta\in\pi_1(E_h,*)$ is proved analogously. We divide the line integrals in $\tht^\pm_{\gamma\delta}$ and the surface integral in $\phi_{\gamma\delta}$ into parts where $t$ or $s$ lies in $\gamma$ or in $\delta$. We use also $\lambda_{\gamma\delta}=\lambda_{\gamma}\lambda_{\delta}$.\\
\end{proof}

\begin{proposition}
  \label{pr:gabilim}
  Let $\xi(t)=(1-x_1)(t)$  with the initial value $\xi(0)=1-x_1^{(1)}$ (see Figure \ref{fig:dreg}). We have the following integral formula for the generalized Abelian integral
  \begin{equation}
    \label{prgabi}
    \gabi_\gamma (h) = (e^{-a T\gamma} -1)^{-1} \Big(\int_{-T_\gamma/2}^{T_\gamma/2}\xi(t) e^{at}\dr t\Big)^2 + \int_{-T_\gamma/2}^{T_\gamma/2}\dr t \int_{-T_\gamma/2}^{t}\dr s\ \xi(s)\xi(t)e^{a(t-s)}.
  \end{equation}
As $h\to 0^+$ these integrals have finite limits:
\begin{equation}
  \label{intlims}
  \begin{split}
    \int_{-T_\gamma/2}^{T_\gamma/2}\xi(t) e^{at}\dr t \longrightarrow \frac{\pi a}{\sin(\pi a/2\sqrt3)}=\sqrt{C_1},\\
    \int_{-T_\gamma/2}^{T_\gamma/2}\dr t \int_{-T_\gamma/2}^{t}\dr s\ \xi(s)\xi(t)e^{a(t-s)} \longrightarrow C_0,
  \end{split}
\end{equation}
where $C_0,C_1$ are as defined in Proposition \ref{pr:genabi}.
\end{proposition}
\begin{proof}
The value of generalized Abelian integral $\gabi_\gamma$ does not depend on the ``shift'' of parametrization (e.g. $t\mapsto t+\tfrac{T_\gamma}2$), but values of integrals $\phi_\gamma$ and $\tht^\pm_\gamma$ depend. We choose Hamiltonian time parameter in such a way that $x_1(0)=x_1^{(2)}$, where $x^{(1)}_h<x^{(2)}_h<x^{(3)}_h$ are real roots of the polynomial $x_1^3-3 x_1 +2=h$ (see Figure \ref{fig:dreg}). Thus 
\[
(1-x_1)(t + T_\gamma/2)=\xi(t)
\]
and so, using formula (\ref{gabirexplicit}) and formulas (\ref{thtpdef}), we get formula (\ref{prgabi}). 

To determine the asymptotic expansion we notice that the singular curve $E_0=\{H(x)=0\}=\{x_2^2=(x_1-1)^2(x_1+2)\}$ is rational. The Hamiltonian parametrization of the limit loop $\gamma_0$ can be explicitly calculated:
\begin{equation}
  \label{xi}
    \xi(t) \longrightarrow \xi_0(t)=\frac3{\cosh^2(\sqrt3 t)},\qquad -\infty<t<\infty
\end{equation}
as $h\to 0^+$. Recall that $T_\gamma(0)=\infty$. Substituting these values to integrals in (\ref{prgabi}) we get the following limits
\[
  \int_{-T_\gamma/2}^{T_\gamma/2}\xi(t) e^{at}\dr t \longrightarrow \int_{-\infty}^{\infty} \frac3{\cosh^2(\sqrt3 t)} e^{it(a/i)}\dr t = \sqrt{2\pi} \mathcal{F}(\xi_0)(a/i)
\]
where $\mathcal{F}$ denotes the Fourier transform. Since $\mathcal{F}(\tfrac1{\cosh^2})(k)=\sqrt{\tfrac{\pi}{2}}\tfrac{k}{\sinh(k\pi/2)}$ (see \cite{gr}, Integral 3.982.1 for example) we find the value
\[
\frac{i \pi a}{\sinh(i\pi a/2\sqrt3)}=\frac{\pi a}{\sin(\pi a/2\sqrt3)}=\sqrt{C_1}.
\]

To determine the limit of the second integral
\[
  \int_{-T_\gamma/2}^{T_\gamma/2}\dr t \int_{-T_\gamma/2}^{t}\dr s\ \xi(t)\xi(s)e^{a(t-s)} \longrightarrow \int_{-\infty}^{\infty}\dr t \int_{-\infty}^{t}\dr s\ \xi_0(t)\xi_0(s)e^{a(t-s)} 
\]
we substitute $u=t-s$, use the symmetry $\xi_0(-t)=\xi_0(t)$ and the Parsival identity; then we get
\begin{align*}
\int_{-\infty}^{\infty}\dr t \int_0^{\infty}\dr u\ \xi_0(t)\xi_0(u-t)e^{a(u)} &= \int_{-\infty}^{\infty}\dr k \Big(\mathcal{F}(\xi_0)(k)\Big)^2 \mathcal{F}\Big(e^{au}\chi_{[0,\infty)} (u)\Big)(k)=\\
&=\frac{3\sqrt{2}i}{\pi^{3/2}}\int_{-\infty}^{\infty} \frac{k^2}{\sinh^2 k}\frac{\dr k}{k-\pi i (a/2\sqrt3)}.
\end{align*}
To evaluate the latter integral, which has the form
\begin{equation}
  \label{fint}
  F(w) = \int_{-\infty}^{\infty} \frac{z^2}{\sinh^2 z}\,\frac{\dr z}{z-\pi i w}, \qquad \real w < 0,
\end{equation}
we must use the logarithmic derivative of the Euler $\Gamma$ function, i.e. the function $\polg=(\log\Gamma)'=\tfrac{\Gamma'}{\Gamma}$. 

Integrating by parts we obtain
\begin{align*}
  F(w) &= \lim_{R\to+\infty}\int_{-R}^{-R^{-1}}\hspace{-1.5em}+\int_{R^{-1}}^{R} (\sgn z-\coth z)'\, \frac{z^2}{z-\pi iw}\dr z = \\
  &=-2\pi i w +\pi^2 w^2\; \lim_{R\to+\infty}\int_{-R}^{-R^{-1}}\hspace{-1.5em}+\int_{R^{-1}}^{R} \frac{\coth z}{(z-\pi i w)^2}\dr z.
\end{align*}
Next, we integrate the function $\tfrac{\coth z}{(z-\pi i w)^2}$ along the contour consisting of the segment $[-R,-R^{-1}]$ followed by the semicircle $R^{-1}e^{i\varphi},\; \varphi\in[\pi,2\pi]$ and segments: $[R^{-1},R]$, $[R,R+i(N+\tfrac12)\pi]$, $[R+i(N+\tfrac12)\pi,-R+i(N+\tfrac12)\pi]$, $[-R,-R+i(N+\tfrac12)\pi]$, where $N\in\bbN$. Using the residue formula and passing to the limit $R,N\to\infty$ we deduce
\newcommand{\res}{\mathrm{Res}}
\begin{multline*}
  \lim_{R\to+\infty}\int_{-R}^{-R^{-1}}\hspace{-1.5em}+\int_{R^{-1}}^{R} \frac{\coth z}{(z-\pi i w)^2}\dr z + \frac{\pi i}{(\pi iw)^2} = 2\pi i\sum_{n=0}^\infty \res_{i\pi n}\Big(\frac{\coth z}{(z-\pi i w)^2}\Big)=\\
  = -\frac{2i}{\pi}\,\sum_{n=0}^\infty \frac1{(n-w)^2} = -\frac{2i}{\pi}\,\polg(-w);
\end{multline*}
in the identification of the latter sum we used \cite{gr}, formula 8.363.8. Finally we have
\[
 F(w) = \pi i \Big(1-2 w-2w^2\polg'(-w)\Big)\qquad \text{for}\quad \real(w)<0,
\]
and so the second of the limits (\ref{intlims}) follows.\\
\end{proof}

Now we investigate the monodromy properties of the generalized Abelian integral $\gabi_\gamma$. We shall need the following
\begin{lemma}[\cite{bo}]
  \label{lem:algid}
For $W,W'\in\trg$ we have
\[
\gabif(W\cdot W') = \gabif(W)+ \gabif(W') + \frac{|W|^2\,|W'|^2} {(|W|^2 -1)(|W'|^2 -1)(|W|^2\,|W'|^2 -1)} \tlgabif([W,W']),
\]
where $[W,W']=W\,W'\,W^{-1}\,(W')^{-1}$ is the commutant and 
\[
\tlgabif(W) = (|W|^2-1)\;\gabif(W);
\]
(for $|W|=1$ we have $\tlgabif(W)=\tht^+\,\tht^-$ in terms of (\ref{gform})).
\end{lemma}
\begin{proof}
The proof relies on direct calculations.\\
\end{proof}

\begin{corollary}
\label{cor:psian}
  The function $\gabi_\gamma(h)$ near $h=0$ has the following form
  \begin{equation}
    \label{psian}
    \gabi_\gamma(h) = \varphi_1(h) + \tfrac1{2\pi i}\gabi_\delta(h) \log h - \frac{\lambda_\delta^2(h)}{(\lambda_\delta^2(h) - 1)^2}\tlgabi_{[\gamma,\delta]}(h)\cdot \frac1{\lambda_\gamma(h) - 1},
  \end{equation}
where $\delta$ is the second cycle in $\pi_1(E_h,*)$ (see the proof of Lemma \ref{lem:abi}), $\tlgabi_{[\gamma,\delta]}(h) = \tlgabif(\rho([\gamma,\delta])$. The functions $\varphi_1$, $\gabi_\delta$ and $\tlgabi_{[\gamma,\delta]}$ are holomorphic near $h=0$.
\end{corollary}
\begin{remark}
\label{rk:tlgabi}
One can prove that the function $\tlgabi_{[\gamma,\delta]}$ is constant
\[
\tlgabi_{[\gamma,\delta]}(h) = (2\pi a)^2.
\]
Indeed, since the contour $[\gamma,\delta]$ is monodromy invariant (see proof of Corollary \ref{cor:psian} below), the function $\tlgabi_{[\gamma,\delta]}$ is meromorphic on whole $\bbC$ with possible poles in $h=0,4$. We know, by Proposition \ref{pr:gabilim}, that it is bounded as $h\to 0$. Similarly one shows that it is bounded as $h\to 4$. These calculations are analogous to proof of the first limit in (\ref{intlims}). We also check that $\tlgabi_{[\gamma,\delta]}$ is bounded as $h\to\infty$. Thus this function has to be constant; its value we calculate by passing to the limit $h\to0$ and comparing respective terms in (\ref{prgabi}) and (\ref{psian}).
\end{remark}

\begin{proof}[Proof of Corollary \ref{cor:psian}]
  The Picard-Lefschetz formula (\ref{piclef}), Theorem \ref{th:bo} and Lem\-ma~\ref{lem:algid} imply that
  \begin{equation}
    \label{mongabi}
    \mon_0\gabi_\gamma =\gabif(\rho(\gamma)\cdot\rho(\delta)) = \gabi_\gamma+\gabi_\delta+\frac{\lambda_\gamma^2\lambda_\delta^2}{(\lambda_\gamma^2 - 1)\,(\lambda_\delta^2 - 1)\,(\lambda_\gamma^2\lambda_\delta^2 - 1)} \tlgabi_{[\gamma,\delta]}
  \end{equation}
and $\mon_0\gabi_\delta=\gabi_\delta$. Next the following monodromy relations follows the Picard-Lefschetz formula (\ref{piclef})
\begin{gather*}
  \mon_0 \lambda_\gamma = \lambda_\gamma\, \lambda_\delta,\qquad \mon\lambda_\delta=\lambda_\delta,\\
  \mon_0[\gamma,\delta] = (\gamma\delta)\cdot\delta\cdot(\gamma\delta)^{-1}\cdot\delta^{-1}=[\gamma,\delta].
\end{gather*}
Therefore $\gabi_\delta$, $\tlgabi_{[\gamma,\delta]}$ and $\lambda_\delta$ are locally single-valued functions of $h$. Since they are bounded (see Proposition \ref{pr:gabilim}), they must be holomorphic. Now 
\begin{align*}
  \mon_0 \Big(\tfrac1{2\pi i} \gabi_\delta\, \log h \Big) &= \tfrac1{2\pi i} \gabi_\delta\, \log h + \gabi_\delta,\\
  \mon_0 \Big(-\frac{\lambda_\delta^2\, \tlgabi_{[\gamma,\delta]} }{(\lambda_\delta^2 - 1)^2\, (\lambda_\gamma^2 - 1)}\Big) &= -\frac{\lambda_\delta^2\, \tlgabi_{[\gamma,\delta]} }{(\lambda_\delta^2 - 1)^2\, (\lambda_\gamma^2\lambda_\delta^2 - 1)} =\\
  = \Big(-\frac{\lambda_\delta^2\, \tlgabi_{[\gamma,\delta]} }{(\lambda_\delta^2 - 1)^2\, (\lambda_\gamma^2 - 1)}\Big) & + \frac{\lambda_\gamma^2\lambda_\delta^2\, \tlgabi_{[\gamma,\delta]}}{(\lambda_\gamma^2 - 1)\,(\lambda_\delta^2 - 1)\,(\lambda_\gamma^2\lambda_\delta^2 - 1)}.
\end{align*}
Therefore the function $\varphi_1$ defined by (\ref{psian}) is single-valued. Since $(\lambda_\gamma^2 -1)$ and $(\lambda_\delta^2 -1)$ are separated from zero and the function $\gabi_\gamma$ is bounded (see Proposition \ref{pr:gabilim}), the function $\varphi_1$ is holomorphic.\\
\end{proof}

Corollary \ref{cor:psian} allows to finish the proof of Proposition \ref{pr:genabi}. The holomorphic function $\varphi_1$ is defined in Corollary and $\varphi_2$, $\varphi_3$ can be read from (\ref{psian}):
\begin{align*}
  \varphi_2 &= \tfrac1{2\pi i} \gabi_\delta,\\
  \varphi_3 &= -\frac{\lambda_\delta^2\, \tlgabi_{[\gamma,\delta]} }{(\lambda_\delta^2 - 1)^2}=-\Big(\frac{2\pi a \lambda_\delta }{\lambda_\delta^2 - 1}\Big)^2
\end{align*}
(the latter equality follows from Remark \ref{rk:tlgabi}).

Since $T_\gamma = -\tfrac1{2\sqrt3} \log h + O(1)$, we have the leading term of expansion
\[
(\lambda_\gamma^2 - 1)^{-1}=(e^{-a T_\gamma} - 1)^{-1} = h^{-a/2\sqrt3}+ ...
\]
\qed
\subsection{Estimates}
\label{sec:estim}

In this subsection we show that the zeroes $h_n$ of the generalized Abelian integral (see Corollary \ref{cor:zeroes}) generate corresponding limit cycles of the system (\ref{system}), provided $\ve$ is sufficiently small. 

Recall that the problem of limit cycles of (\ref{system}) is reduced to the problem of limit cycles of the following planar system
\begin{equation}
  \label{plsystem}
  \dot x = X_H(x) + \ve\,\real\Big(\overline{\bz}\,G(x,\ve)\Big) e_2,
\end{equation}
where the function $G(x,\ve)$ is defined via the invariant surface $L_\ve = \{y=\ve G(x,\ve)\}$, which is a graph of $G(\cdot,\ve)$.

At the moment we do not even know whether the invariant surface exists. Indeed, the normal hyperbolicity conditions are \emph{not} satisfied: the eigenvalues in the normal direction are $\lambda_{3,4} = -\sqrt3 \pm i \sqrt3$, whereas the eigenvalues at the saddle point $x=(1,0),y=0$ in the $x$-direction are $\pm 2\sqrt3$ (compare \cite{bozoell,bozo3d,hps,ni}). We should do two things:
\begin{enumerate}
\item prove the existence of the invariant surface, 
\item estimate the discrepancy $G(x,\ve)- H^4\, g(x)$, where $g(x)$ is the solution to the normal variation equation given in (\ref{nveq}).
\end{enumerate}
In both tasks the crucial role is played by the following Lemma. Let us recall the notation related to the elliptic Hamiltonian $H(x)=x_1^3 -3x_1 - x_2^2+2$. The basin $D\subset\bbR^2$ (see Figure \ref{fig:dreg}) is filled with closed orbits of the Hamiltonian vector field $X_H$. 
\begin{lemma}
\label{lem:c1inv}
Let $U\supset D\times \{0\}$ be an open neighborhood in $\bbR^2\times \bbC$ and $V_\ve$ be the following vector field in $U$
\begin{equation}
  \label{vve}
V_\ve \left\{  \begin{aligned}
    \dot x &= X_H + H^k\,\real (\overline{\bz}\,y) v_o + Q(x,y;\ve)\\
    \dot y &= a y + B(x,y;\ve),
  \end{aligned} \right.
\end{equation}
where $k\geq 0$, $a=-\sqrt3+i\sqrt3$, $\bz\in\bbC$, $v_0\in\bbR^2$ and $Q,B$ are functions of class $C^2(U)$ satisfying the following
\[
\begin{split}
  |Q| &\leq \const\cdot  \ve|H(x)|^2,\\
  |B| &\leq \const\cdot  \ve|H(x)|^2.
\end{split}
\]
Then, for sufficiently small $\ve$ there exists a unique invariant surface of $V_\ve$:
\[
L_\ve = \{(x,y):\quad x\in D, \ y=\ve\, G(x,\ve)\}.
\]
The function $G(\cdot,\ve)$ prolongs by zero to a $C^1$ function on a neighborhood of $D$ in $\bbR_x^2$.
\end{lemma}

\begin{proof}
We shall prove that the Poincar\'e return map associated to vector field $V_\ve$ satisfies the normal hyperbolicity condition. 

In a neighborhood of the center critical point $(x=(-1,0),y=0)$ for unperturbed vector field $V_0$, the system is normally hyperbolic and so the invariant surface exists. In the further proof we shall concentrate on the neighborhood of separatrix $\gamma_0=\{(x,y):y=0,\quad H(x)=0\}$ (see Figure \ref{fig:estimaux}).
\begin{figure}[htbp]
  \centering
  \input{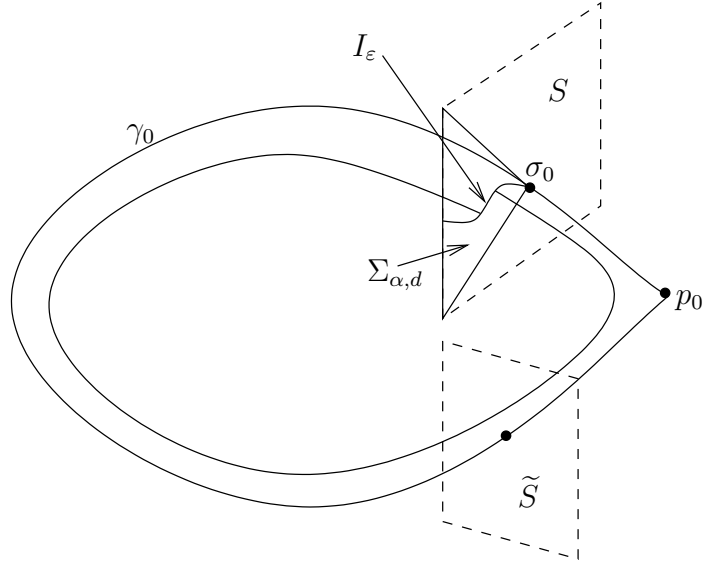}
  \caption{The Poincar\'e return map $\poi_\ve$ defined via trajectories of $V_\ve$ and its invariant unit $I_\ve$.}
  \label{fig:estimaux}
\end{figure}

The non-degenerate critical point $p_0=((1,0),0)$, located on this separatrix, is preserved after the perturbation. Let us choose a 3-dimensional hypersurface $S$ transversal to $\gamma_0$ and close to the singular point $p_0$. Let
\begin{align*}
    \Sigma_{\alpha,d} &= S\cap \{(x,y):\quad |y|^2 \leq \alpha^2 (H(x))^2, \quad H(x) < d\}\\
\intertext{be a sector in $S$ with vertex}
    \sigma_0 &= \Sigma_{\alpha,d}\cap \gamma_0 = \{(x,0)\in\Sigma_{\alpha,d}:\quad  H(x)=0\}.
\end{align*}
\begin{lemma}
\label{lem:poinc1}
For sufficiently small $\alpha >0,\ d_2>d_1>0$ and $\ve\in(-\ve_0,\ve_0)$, the Poincar\'e return map
\begin{equation}
  \label{pvempa}
  \poi_\ve : \Sigma_{\alpha,d_1} \rightarrow \Sigma_{\alpha,d_2}
\end{equation}
is a diffeomorphism onto the image and prolongs to a map of class $C^1$ in point $\sigma_0$.
\end{lemma}

Now we finish the proof of Lemma \ref{lem:c1inv}. The unperturbed Poincare map $\poi_0$ is the identity on the invariant segment $I_0=\Sigma_{\alpha,d_1}\cap\{y=0\}$. Thus $\poi_0$ is normally hyperbolic on $I_0$, since we have strong contraction in the normal direction. In virtue of the Hirsh Pugh Shub Theorem \cite{hps}, for sufficiently small $\ve$ there exist the unique invariant embedded interval $I_\ve$ close to $I_0$; it is of class $C^1$. Considering Hamiltonian as the parameter on $I_0$, we get 
\[
I_\ve =\{y=\ve F(h,\ve), \ h\in [0,\delta)\},\quad F\in C^1([0,\delta)\times (-\ve_0,\ve_0)).
\]
The surface $S_\ve$ spanned by trajectories of $V_\ve$ passing through $I_\ve$ is $V_\ve$-invariant, due to the invariance of $I_\ve$ under the Poincar\'e map $\poi_\ve$. The form of invariant interval $I_\ve$ and the form of vector field $V_\ve$ implies that the surface $S_\ve$ is graph of a $C^1(D)$ function:
\[
L_\ve=\{(x,y):x\in D, \quad y=\ve G(x,\ve)\}.
\]
We prove that $G$ can be extended by zero outside $D$.

We check that, providing the assumptions of Lemma \ref{lem:c1inv} hold, the set
\[
\{(x,y):x\in D, |y|^2 \leq R^2 \ve^2 |H(x)|^4\},
\]
for $R$ big enough, is invariant for $V_\ve$. Namely, denoting all constants by $C$, we calculate
\begin{multline}
\label{invcal}
  V_\ve(|y|^2 - R^2 \ve^2 |H|^4) |_{|y| =R\, \ve\, |H|^2} =\\
= 2 \real\Big[\overline{y}\; (a\, y+B)\Big] -4 R^2 \ve^2 H^3\Big(H^3\real(\overline{\bz}\,y\, <\dr H,v_0> + <\dr H,Q>)\Big) \\
  \leq 2R^2\,\ve^2\, H^4 \Big(-\rho + \tfrac{C}{R} + 2R\,\ve\,C\,|H|^{1+k}+ \ve\,C\,|H|\Big).
\end{multline}
Since the latter expression is $\leq 0$ (for sufficiently large $R$), the considered subset is $V_\ve$-invariant.

Thus 
\[
|G(x,\ve)|\leq \const \cdot |H(x)|^2
\]
and the function $G(x,\ve)$ can be prolonged by zero to a function of class $C^1$.\\
\end{proof}

\begin{proof}[Proof of Lemma \ref{lem:poinc1}]
By calculations analogous to (\ref{invcal}) we find that
\[
  V_\ve(|y|^2 - \alpha^2 H^2) |_{|y| =\alpha\, |H|} \leq 2\alpha^2\, H^2 \Big(-\rho + \tfrac{\ve\,|H|}{\alpha} + \alpha\,C\,|H|^k + C\ve\,|H|\Big).
\]
Thus for sufficiently small $\alpha$ and $\ve$ the subset
\[
\{(x,y):x\in D,\quad |y|^2\leq \alpha^2 (H(x))^2\}
\]
is $V_\ve$-invariant. Moreover, the separatrix $\gamma_0$ is also $V_\ve$-invariant. This proves, that the Poincar\'e return map defines the dyffeomorphism (\ref{pvempa}) which is of class $C^1$ outside the border.

Thus it remains to show that $\poi_\ve$ can be prolonged to the $C^1$ map in the point $\sigma_0$. We choose an additional, auxiliary, 3-dimensional hypersurface $\widetilde{S}$ transversal to $\gamma_0$, close to $p_0$, which lies ``on another side'' with respect to the point $p_0$ (see Figure \ref{fig:estimaux}). The return map $\poi_\ve$ is the composition $\poi_\ve = \poi^s_\ve \circ \poi^r_\ve$ of the correspondence maps
\[
\poi^s_\ve : S\rightarrow \widetilde{S},  \qquad \text{and}\qquad \poi^r_\ve : \widetilde{S}\rightarrow S,
\]
defined by trajectories near the singular point $p_0$ and trajectories near the regular part of $\gamma_0$ respectively. The regular map naturally extends to the $C^1$ map (even $C^2$ in fact) as the flow of the non-vanishing, $C^2$ vector field $V_\ve$. To analyze the singular part $\poi^s_\ve$, we use the following theorem of H. Belitskii.
\begin{theorem}[\cite{bel}]
Let $\Lambda\in \mathrm{End}(\bbR^n)$  be a linear endomorphism whose eigenvalues $(\lambda_1,\ldots,\lambda_n)$ satisfy
\[
\mathrm{Re}\lambda_i \neq \mathrm{Re}\lambda_j + \mathrm{Re}\lambda_k
\]
for all $i$ and $j,k$ such that $\mathrm{Re}\lambda_j \leq 0\leq \mathrm{Re}\lambda_k$. Then any $C^2$ differential system
\[
\frac{\dr x}{\dr t} = \Lambda x + f(x), \qquad f(0)=0=f'(0)
\]
is in the neighborhood of $0$ $C^1$-equivalent to the linearization.
\end{theorem}
The eigenvalues of the linearization of our vector field $V_\ve$ in $p_0$ are $\pm2\sqrt3,\, -\sqrt3\pm i\sqrt3$, so $V_\ve$ satisfies the assumptions of the Belitskii theorem. In suitable coordinates $(u,v)$, associated with the linearization of $V_\ve$ in the neighborhood of $p_0$, the correspondence map $\poi^s_\ve$ has the form 
\[
\poi^s_\ve(u,v)=(u,C\  u^\beta\, v),\quad u\in\bbR_+,\quad v\in (\bbC,0),\quad \beta= \tfrac12-\tfrac{i}2.
\]
The restriction to $\Sigma_{\alpha,d}$ corresponds to the restriction to the set $\{|v|^2\leq \widetilde{\alpha}(u)\, u$, $u\in\bbR_+$, $\widetilde{\alpha}(0)>0\}$. In such region the map $\poi^s_\ve$ is of class $C^1$ in $\sigma_0$; we have $(u,v)(\sigma_0)=(0,0)$ and $(\poi^s_\ve)'(0,0)=\left(\begin{smallmatrix}1&0\\ 0&0\\ \end{smallmatrix}\right)$. Thus, the thesis of the lemma follows.\\
\end{proof}

Now we can prove the existence of the invariant surface and estimate the distance to its linear approximation $H^4g$.
\begin{proposition}
\label{pr:discr}
For sufficiently small $\ve$, there exists the invariant surface 
\begin{equation}
  \label{invsurf}
L_\ve=\{y=\ve\, G(x,\ve),\quad x\in D \}
\end{equation}
of the system (\ref{system}). The function $G$ is of class $C^1$ and the distance to the linearization $(H^4g)(x)=G(x,0)$ is bounded by
 \begin{align}
    |G - H^4g| &\leq C |\ve| |h|^5, \label{Gg} \\
    |G'-(H^4g)'| &\leq C |\ve| |h|^4, \label{Ggprim}
  \end{align}
where $G'=G,_x$ is the derivative with respect to $x$.
\end{proposition}
\begin{proof}
The existence of the invariant surface of class $C^1$ and the form (\ref{invsurf}) is a direct consequence of Lemma \ref{lem:c1inv}. 

To show the bounds (\ref{Gg}), (\ref{Ggprim}), we make a coordinate change
\[
(x,y) \longmapsto (x,z),\qquad z=\Big(y-\ve\, H^4g(x)\Big)/(\ve\, H^4).
\]
The system (\ref{system}) takes the form
\begin{equation}
  \label{systemxz}
  \left\{
    \begin{aligned}
      \dot{x} =& X_H + H^4\,\real(\overline{\bz}\,z)\,\ve\,e_2 + \ve\,H^4\,\real(\overline{\bz}\,g)\,e_2,\\
      \dot{z} =& a\,z + 4\ve\,H^3 (z+g) \real(\overline{\bz}\,(z+g))\, 2x_2 - \ve\,H^4\tfrac{\partial g}{\partial X_2}\,\real(\overline{\bz}\,(z+g)).
    \end{aligned}\right.
\end{equation}
Using the integral formula (\ref{gint}) for the function $g$ we deduce that it is bounded, $g\leq C$. Since the function $g$ prolongs to the (multivalued), holomorphic function ramified along the singular curve $\gamma_0$, the following bounds for derivatives of $g$ hold
\begin{equation}
\label{gbounds}
|g^{(k)}|\leq C_k |H|^{-k}.
\end{equation}
Using this one can check that the system (\ref{systemxz}) satisfies the assumptions of Lemma \ref{lem:c1inv}. Thus, the invariant surface has the form
\[
z=\frac{\ve\,G-\ve\,H^4g}{\ve\,H^4}=\ve\,U(x,\ve)
\]
and the function $U$ prolongs by zero to a $C^1$ function on a neighborhood of $D$. Thus, the function $U$ satisfies the estimates
\[
|U|\leq C|H|, \qquad |U'|\leq C,
\]
which are equivalent to (\ref{Gg}) and (\ref{Ggprim}).\\
\end{proof}

Now we show that the generalized Abelian integral $J(h)$ is a good approximation of the Poincar\'e return map and so the zeroes of $J(h)$ generate limit cycles for sufficiently small $\ve$.
\begin{proposition}
\label{pr:jlead}
Let $\Delta H (h,\ve)$ be the increment of the Hamiltonian after the first return of the system $V_\ve$ restricted to the invariant surface $L_\ve$. Then, there exists a constant $C$ such that
\begin{align}
  |\Delta H -\ve J| &\leq C\, \ve\, |h|^5, \label{dist}\\
  |\partial_h (\Delta H) -\ve J'(h)| &\leq C\,\ve\, |h|^4\ |\log h|. \label{distprim}
\end{align}
\end{proposition}

\begin{proof}
Here we study the phase curves of the 2-dimensional vector field (\ref{plsystem}) i.e.
\[
W_\ve \colon = V_\ve|_{L_\ve} = X_H + \ve Re(\overline{\bz}\, G)\, \partial_{x_2}.
\]

We fix the segment $I=\{(x_1,0):x_1\in [-2,-1)\}$ transversal to the Hamiltonian flow. We denote by $\ircv_\ve(t,h)$ the integral curves of $W_\ve$ which start and finish at $I$. They satisfy
\begin{equation}
\label{auxbounds}
  \begin{split}
    \dot{\ircv}_\ve (t,h) &= X_H + \ve\,\real(\overline{\bz}\,G)\; e_2,\\
    \ircv_\ve (0,h) &\in I,\qquad \ircv_\ve (T_\ve,h)\in I  \\
    H(\ircv_\ve(0,h)) &= h.
  \end{split}
\end{equation}
For $\ve=0$ the curve $\ircv_0$ is the oval $\{H=h\}$ and for $\ve$ non-zero but small it is a small perturbation of $\ircv_0$:
\[
\ircv_\ve (t,h) = \ircv_0 (t,h) + \ve \dis(t,h;\ve).
\]
Above $T_\ve=T_\ve(h)$ is the time of the first return to the unit $I$.
\begin{lemma}
  \label{lem:bound}
  There exist a constant $C$ and the positive, small constant $\nu $ such that the following estimates hold:
  \begin{align}
    |\dis| &\leq C |h|^4\; e^{(2\sqrt3 +\nu) t}, \label{bndF}\\
    |\partial_h \dis| &\leq C |h|^3 \ e^{(2\sqrt3 +\nu) t},\label{bndFprim}\\
    |T_\ve - T_0| &\leq C\,\ve\, |h|^{3}.\label{bndtt}
  \end{align}
\end{lemma}
\begin{proof}
We use the scalar product $x\cdot x' = 3 x_1\,x_1' + x_2\,x_2'$ and $|x|=\sqrt{x\cdot x}$.

It follows from the equation (\ref{auxbounds}) that the function $\dis$ satisfies the following initial value problem
\begin{equation}
\label{auxF}
  \left\{
    \begin{aligned}
      \dot{\dis} &= \widetilde{\dr X}_H\, \dis + \real(\overline{\bz}\,G(\ircv_0+\ve\,\dis))\, e_2,\\
      \dis(0,h;\ve) &= 0,
    \end{aligned}\right.
\end{equation}
where $\widetilde{\dr X}_H\, \dis = \tfrac1\ve ((X_H(\ircv_0+\ve\,\dis) - X_H(\ircv_0))$. We have $\widetilde{\dr X}_H = \dr X_H (\ircv_0+\theta\,\ve\,\dis)$, for some $\theta\in(0,1)$. Hence
\begin{equation*}
  \widetilde{\dr X}_H = \left(
    \begin{smallmatrix}
      0& -2\\ -6x_1&0
    \end{smallmatrix}\right),\qquad x_1\in [-2,1].
\end{equation*}
 Moreover, using estimates (\ref{Gg}), (\ref{Ggprim}), (\ref{gbounds}) we get $|G(\ircv_0+\ve\,\dis)|\leq C_1\,h^4$. For the solution $\dis$ to the equation (\ref{auxF}) we have
\begin{multline*}
\tfrac{\dr}{\dr t} | \dis|^2 = 2 | \dis| \; \tfrac{\dr}{\dr t} |\dis | = 2 \dis\cdot\widetilde{\dr X}_H\, \dis + 2 \dis\cdot\real(\overline{\bz}\,G(\ircv_0+\ve\,\dis))\, e_2 =\\
= -12(x_1+1)\, \dis_1\,\dis_2 + 2 \dis_2\, \real(\overline{\bz}\,G(\ircv_0+\ve\,\dis)) \leq 4\sqrt3 | \dis |^2 + 2 C_2 h^4 | \dis|.
\end{multline*}
Therefore $\tfrac{\dr}{\dr t} \left(| \dis|\right) \leq 2\sqrt3 | \dis| + C_2 h^4,\quad |\dis| (0) = 0$ and the Gronwall inequality \cite{har} gives the bound (\ref{bndF}).

Since the difference of flows $\ve\,\dis$ after the Hamiltonian period $T_0$ is $\leq \widetilde{C}|h|^{3}$ and the ``velocity'' $|V_\ve|\sim 1$, the difference of periods $|T_\ve-T_0|$ satisfies (\ref{bndtt}).

The derivative $\tfrac{\partial \dis}{\partial_h}$ satisfies the respective linear variation equation related to (\ref{auxF})
\[
\tfrac{\dr}{\dr t}\Big(\tfrac{\partial \dis}{\partial_h}\Big) = \Big(\widetilde{\dr X}_H+\ve\,(\ldots)\Big)\,\tfrac{\partial \dis}{\partial_h} + h^3\,(\ldots)\,(\tfrac{\partial \ircv_0}{\partial_h}),\qquad \tfrac{\partial \dis}{\partial_h}(0,h;\ve)=0,
\]
where we denoted by $(\ldots)$ the bounded terms (see estimations (\ref{bndF},\ref{bndtt},\ref{Ggprim},\ref{gbounds})). Since the flow variation $\tfrac{\partial \ircv_0}{\partial_h}$ of the Hamiltonian field satisfies $\left|\tfrac{\partial \ircv_0}{\partial_h}\right|\leq C\,e^{2\sqrt3s}$, the inequality (\ref{bndFprim}) holds. This finishes the proof of Lemma \ref{lem:bound}.\\
\end{proof}

\noindent
\emph{We continue the proof of Proposition \ref{pr:jlead}.}\\
We split the difference between the Poincar\'e map and the linearization $\ve\, J(h)$ in two integrals $R_1(h,\ve)$ and $R_2(h,\ve)$:
\begin{align*}
  R_1 &= \int_{\gamma_\ve} \real(\overline{\bz}\,(G-H^4g))\; \dr x_1,\\
  R_2 &= \int_{\gamma_\ve} \real(\overline{\bz}\,H^4g)\; \dr x_1 - \int_{\gamma_0} \real(\overline{\bz}\,H^4g)\; \dr x_1.
\end{align*}
We shall show that the estimations (\ref{dist}) and (\ref{distprim}) hold for both $R_1$ and $R_2$.

The inequality (\ref{dist}) is a direct consequence of the bound (\ref{Gg}). 
The difference of $R_1$ in close values of $h$ takes the form
\[
R_1(h+\delta)- R_1(h)  = \int_{\gamma_\ve(h,h+\delta)} \partial_{x_2}\real(\overline{\bz}\,(G-H^4g))\;\dr x_2\wedge\dr x_1 + O(\ve |h-2|^5),
\]
where the integral is taken along the strip $\gamma_\ve(h,h+\delta)$ between $\gamma_\ve(h)$ and $\gamma_\ve(h+\delta)$. The area of this strip is of the same order as the area of the domain $\{X:\ h<H<h+\delta\}$, i.e.\ $\sim \delta\cdot I_0\sim C\,\delta\,|\log h|$ (see proof of Lemma \ref{lem:abi}). So the estimate (\ref{distprim}) follows (\ref{Ggprim}).

To prove the estimate for $R_2$ we use (\ref{auxbounds}) and (\ref{gbounds}):
\begin{multline*}
  |R_2|\leq C_1 \int_0^{T_0}\real(\overline{\bz}\,H^4g) (\ircv_0+ \ve \dis) - \real(\overline{\bz}\,H^4g) (\ircv_0) + \int_{T_0}^{T_\ve} C_2\,\ve\,|h|^4 \leq \\
  \leq C_3\,\ve\, |h|^7 \int_0^{T_0} e^{(2\sqrt3 +\nu)t}\;\dr t + C_4\,|h|^4\,|T_\ve-T_0|\leq C_5\,\ve\, |h|^{6-2\nu} \leq C \ve |h|^5.
\end{multline*}

Similarly, using the following formula for differential of integral
\[
\frac{\partial}{\partial h} \int_{\ircv_\ve} \omega = \int_{\ircv_\ve} i_{\partial \ircv_\ve/\partial_h}\;\dr \omega + \omega \left(\left.\tfrac{\partial\ircv_\ve}{\partial_h}\right|_{t=0}\right) -  \omega \left(\left.\tfrac{\partial\ircv_\ve}{\partial_h}\right|_{t=T_\ve}\right)
\]
and the bounds (\ref{auxbounds}), (\ref{gbounds}), we get
\[
|R_2|\leq C\,\ve\, |h|^5,\qquad |\partial_h R_2|\leq C\,\ve\, |h|^4.
\]
\end{proof}

Now we can finish the proof of the Main Theorem. The restriction of system (\ref{system}) to its invariant surface $L_\ve=\{y=\ve\,G(x,\ve)\}$ has the form (\ref{plsystem}). The increment $\Delta H=\poi(h)-h$, associated with the Poincar\'e map $\poi(h)$ (on a section transversal to $\gamma_0$), equals (see Proposition \ref{pr:jlead})
\begin{align*}
  \Delta H(h) &= \ve J(h) + O(|\ve|\,|h|^5),\\
  (\Delta H)'(h) &= \ve J'(h) + o(|\ve|\,|h|^{4-1/4}).
\intertext{Since we also have (see proof of Corollary \ref{cor:zeroes})}
  J(h) &= R h^{4+1/2} \cos(\log\sqrt h -\alpha_0) + o(h^{4-1/4}),
\end{align*}
any simple zero of $J(h)$, which is sufficiently close to $0$ generates, by the Implicit Function Theorem, a simple zero of $\Delta H$. Thus the sequence $h_n\to 0^+$ of simple zeroes of $J(h)$ (see Corollary \ref{cor:zeroes}) guarantees the existence of \emph{infinite sequence} $\widetilde{h}_n\to 0^+$, $n\geq N_0$ of simple zeroes of the increment $\Delta H$. Any simple zero of $\Delta H$ corresponds to a limit cycle of (\ref{system}).

The proof of Main Theorem is now complete.\\
\qed

\end{document}